\title{On winning fast in Avoider-Enforcer games}
\author{{J\'anos Bar\'at\thanks{Department of Computer Science and Systems Technology, University of Pannonia, Egyetem u. 10, 8200 Veszpr\'em, Hungary.
Research is supported by OTKA Grant PD~75837.}}
\quad and \quad {Milo\v{s} Stojakovi\'{c} \thanks{Department of Mathematics and
Informatics, University of Novi Sad, Serbia; Partly supported by
Ministry of Science and Environmental Protection, Republic of
Serbia, and Provincial Secretariat for Science, Province of
Vojvodina.}}}

\documentclass[12pt]{article}
\usepackage{amsmath,amssymb,latexsym,a4}
\usepackage[dvipdfm]{graphicx}

\newtheorem{theorem}{Theorem}[section]

\newtheorem{question}[theorem]{Question}

\DeclareMathOperator{\littleo}{o}
\DeclareMathOperator{\BigO}{O}

\begin{document}
\maketitle

\begin{abstract}
We analyze the duration of the unbiased Avoider-Enforcer game for three basic positional
games. All the games are played on the edges of the complete graph on $n$ vertices, and
Avoider's goal is to keep his graph outerplanar, diamond-free and $k$-degenerate,
respectively. It is clear that all three games are Enforcer's wins, and our main interest lies
in determining the largest number of moves Avoider can play before losing.

Extremal graph theory offers a general upper bound for the number of Avoider's moves. As it
turns out, for all three games we manage to obtain a lower bound that is just an additive
constant away from that upper bound. In particular, we exhibit a strategy for Avoider to keep
his graph outerplanar for at least $2n-8$ moves, being just $6$ short of the maximum possible.
A diamond-free graph can have at most $d(n)=\lceil\frac{3n-5}{2}\rceil$ edges, and we prove
that Avoider can play for at least $d(n)-3$ moves. Finally, if $k$ is small compared to $n$,
we show that Avoider can keep his graph $k$-degenerate for as many as $e(n)$ moves, where
$e(n)$ is the maximum number of edges a $k$-degenerate graph can have.
\end{abstract}

\section{Introduction} \label{sec:intro}

In this paper, we deal with Avoider-Enforcer positional games. For a hypergraph ${\mathcal
F}$, the game ${\mathcal F}$ is played by two players, Avoider and Enforcer. They alternately
claim previously unclaimed vertices of $\mathcal{F}$. Avoider starts, and the game ends when
all vertices have been claimed. Enforcer wins if Avoider has claimed all vertices of some
hyperedge of $\mathcal{F}$. Otherwise Avoider wins. We refer to the vertices of $\cal F$ as
the board, and the hyperedges of $\cal F$ as the losing sets. The recent book~\cite{Beckbook}
by Beck offers a good overview of the topic of positional games. Here, we study games which
are played on the edges of the complete graph on $n$ vertices, that is, the board of
${\mathcal F}$ is always $E(K_n)$.

If we assume that both players play optimally, then each game $\mathcal{F}$ is either an
Avoider's win or an Enforcer's win. A significant part of the previous work done in
combinatorial game theory (see, e.g.,~\cite{winways}) is devoted to the question: Which one of
the two players wins a particular game? Here, we go one step further and address a different
issue -- our hope is to determine not only the winner of a game, but also {\em how fast} is he
able to win.

For a game $\cal F$, let $\tau_E({\mathcal F})$ be the smallest integer $t$ such that Enforcer
has a strategy to win the game ${\mathcal F}$ in at most $t$ moves. We say that $\tau_E({\cal
F}) = \infty$, if the game is an Avoider's win.

For an Avoider-Enforcer game, this type of question was first raised only recently, by Hefetz
et al.\ in~\cite{HKSS}, and it was also addressed in~\cite{AJSS}. On the other hand, an
analogue question for Maker-Breaker games, the more studied Avoider-Enforcer games'
counterpart, has been a topic for some time. We mention here the work of Beck~\cite{Beckfast}
and Peke\v c~\cite{Pekec}, who looked at how fast Maker can win the clique game. Chv\'atal and
Erd\H os~\cite{CE}, and later Hefetz et al.~\cite{HKSS2}, studied the fast winning in
Maker-Breaker Hamiltonicity game.

We would like to emphasize that, generally speaking, results on fast winning in positional
games have an impact on the whole field, as those results can later be used in analysis of
other positional games. Namely, it often happens that an optimal strategy of a player consists
of several stages, and in each of them the player wants to complete a task. In that situation,
a particular task should not only be performed, but performed fast, i.e., in significantly
less moves than the total number of moves at player's disposal.

\subsection{Preliminaries}

All in all, the theory behind Avoider-Enforcer games is less developed than the one behind
Maker-Breaker games. However, when it comes to determining how fast can Enforcer win the game,
a somewhat unexpected help comes from extremal graph theory.

The {\em extremal number} (or Tur\'an number) of a hypergraph
$\cal F$ is defined by $\text{ex}({\cal F}) = \max\left\{|A|:
A\subseteq V({\cal F}),\,\,A\not\in E({\cal F}) \right\}$.
As it was shown in~\cite{HKSS}, if the set of hyperedges of $\cal F$ is a
monotone increasing family of sets, we have
\begin{equation}
\frac{1}{2}\text{\rm ex}({\cal F})+1 \leq  \tau_E({\cal F}) \leq
\text{\rm ex}({\cal F}) +1. \label{e:ob}
\end{equation}
Note that for every game $\cal F$, we can make the set of hyperedges
an increasing family by adding all the supersets of the hyperedges
to the set of hyperedges. This process changes neither the outcome
nor the nature of the game.

Therefore, as soon as we know the extremal number for the game hypergraph, from (\ref{e:ob})
we get the length of the game squeezed between two values which are roughly a factor of two
from each other.

In~\cite{AJSS} and~\cite{HKSS}, the possibilities of Enforcer's fast win for several
well-studied positional graph games were analyzed. As it was shown in~\cite{AJSS}, Avoider can
keep his graph planar for as many as $3n-O(1)$ moves, which is just a constant away from the
upper bound derived from~(\ref{e:ob}). Two other basic positional games are looked at
in~\cite{HKSS}. In the first one, Avoider wants to keep his graph bipartite for as long as
possible, where in the second one his goal is to avoid creating a spanning graph. The duration
of both games is determined quite precisely in both the first and the second order terms. It
turns out that, in both cases, the values are not additive constant away from either of the
bounds in~(\ref{e:ob}).

\subsection{Our results}

In the present paper, we analyze the duration of the Avoider-Enforcer game for three basic
positional games. As we saw in the non-planarity game, in contrast to several other games that
were analyzed, Avoider can keep his graph planar for quite a long time, just constant away
from the upper bound in~(\ref{e:ob}). We were curious as to what are the reasons behind this,
analyzing a fairly similar game -- the game in which Avoider wants to keep his graph
outerplanar for as long as possible. Formally, let ${\cal OP}_n$ be the hypergraph whose
hyperedges are the edge-sets of all non-outerplanar graphs on $n$ vertices. The
relation~(\ref{e:ob}) shows that $n\leq \tau_E({\cal OP}_n)\leq 2n-2$, which leaves $n-1$
possible values for $\tau_E$. We manage to narrow down the choice to just five values.

\begin{theorem} \label{th:op}
$$2n-7\leq \tau_E({\cal OP}_n)\leq 2n-3.$$
\end{theorem}

We see that, similarly to the non-planarity game, the duration of the game is just an additive
constant away from the upper bound obtained from~(\ref{e:ob}). The common feature of the
non-outerplanarity game and the non-planarity game is that in both cases Avoider loses as soon
as his graph contains a certain forbidden minor. Indeed, for outerplanarity these forbidden
minors are $K_4$ and $K_{2,3}$, and for planarity the forbidden minors are $K_5$ and
$K_{3,3}$. We were curious to analyze further the games of this kind. Hence, we turned our
attention to a game where Avoider's goal is to avoid a single forbidden minor in his graph.
The forbidden minor is the diamond, that is, $K_4$ with one edge missing. We note that
\cite{HKSS3} deals with a similar game, where Avoider's goal is to avoid claiming a fixed
minor. However, the game analyzed there is biased, and the main interest is just the final
outcome.

Formally, let ${\cal DF}_n$ be the hypergraph whose hyperedges are the edge-sets of all graphs
on $n$ vertices which contain a diamond minor. As the number of edges in a diamond-free graph
is at most $d(n)=\lceil\frac{3n-5}{2}\rceil$, from~(\ref{e:ob}) we get $\frac{1}{2}d(n)+1 \le
\tau_E({{\cal DF}_n})\le d(n)+1$. In the following theorem, we reduce this interval to three
integers, again an additive constant away from the upper bound.

\begin{theorem} \label{th:df}
$$d(n)-2\le \tau_E({{\cal DF}_n})\le d(n).$$
\end{theorem}

We note that diamond-free graphs are also called cactus graphs, and
it can be shown that they are outerplanar.

A graph $G$ is called $k$-degenerate, if every subgraph of $G$ has a vertex of degree at most
$k$. The degeneracy of a graph is the minimal $k$ such that the graph is $k$-degenerate. Low
degeneracy is a common property of planar and outerplanar graphs; their
degeneracy is at most 5 and 2, respectively. It is known that graph degeneracy plays a key
role in several other positional games on graphs, see, e.g.,~\cite{mSt}.

Here, our aim is to study a game in which Avoider's goal is to keep his graph $k$-degenerate,
for an integer $k$. In a way, it brings all the mentioned games together, as its family of
forbidden graphs, for some values of $k$, contains the aforementioned families of forbidden
graphs.

Formally, let ${\cal D}_n^k$ be the hypergraph whose hyperedges are the edge-sets of all
graphs on $n$ vertices which are not $k$-degenerate. A $k$-degenerate graph with $n$ vertices
can have at most $e(n)=(n-k)k+{k\choose 2}$ edges, and we show that Avoider loses only when he
has claimed more than $e(n)$ edges, assuming that $n$ is large enough compared to $k$.

\begin{theorem} \label{th:kd}
If $k=o(\log n)$, then $\tau_E({\cal D}_n^k)=e(n)+1$.
\end{theorem}

Our graph-theoretic notation is standard and follows that of~\cite{DBW}.
A matching $M$ of a graph $G$ is called {\sl near-perfect} if there are
at most two $M$-unsaturated vertices in $G$.
If $H$ is a graph, we say that a graph $G$ is $H$-free, if $G$
contains no $H$-minor.
Throughout the paper, $\log$ stands for the natural logarithm.

Occasionally, we may work with dynamic sets and notations. For instance, $\cal A$ denotes the
set of edges claimed by Avoider. At the start of a game, it is the empty set. If Avoider
claims the edge $e$ in his $i$-th move, then we change $\cal A$ to be  ${\cal A}\cup e$.


\section{The strategies -- fast winning and slow losing}

\subsection{Keeping the graph outerplanar} \label{sec:outer}

{\it Proof of Theorem~\ref{th:op}.} Assume that Enforcer claims an edge $uv$ in his first
move. If at any time Avoider claims an edge incident to $uv$, say $xv$, then Enforcer claims
the edge $xu$ in the next move. This simple pairing strategy enables Enforcer to prevent
Avoider from claiming any triangle on the edge $uv$. Therefore, Avoider is unable to claim a
maximal outerplanar graph and loses after at most $2n-3$ moves.

Next, we show a strategy for Avoider to keep his graph $\cal A$ outerplanar for $2n-8$ moves.
In his first two moves, Avoider claims two edges of a triangle. We denote the third edge of
this triangle by $m$. Note that ${\cal A} \cup \{m\}$ is a maximal outerplanar graph on three
vertices. For most of the game, Avoider maintains the graph $\cal A$ to consists of a graph
one edge short of a maximal outerplanar graph, and some isolated vertices. He achieves that by
attaching an isolated vertex to the current outerplanar graph in every pair of consecutive
moves.

\begin{figure}[ht]
 \begin{center}
 \includegraphics[scale=0.5]{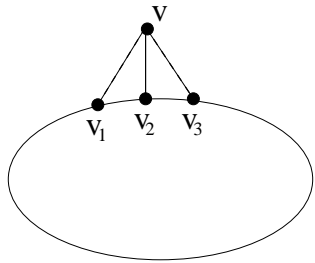}
 \caption{Extension process, which can be performed for good vertices}
\label{f:ext}
\end{center}
\end{figure}

Throughout the game, we denote the outer face of ${\cal A}\cup \{m\}$ by $O_{\cal A}$. An
isolated vertex $v$ in Avoider's graph will be called \emph{bad}, if for every three
consecutive vertices $v_1,v_2,v_3$ on $O_{\cal A}$ at least one of the edges $vv_1$, $vv_2$,
$vv_3$ is claimed by Enforcer. Any other isolated vertex of $\cal A$ is called \emph{good}. A
good vertex can be attached to the current outerplanar graph $\cal A$  in two Avoider's moves.
Namely, if $v_1,v_2,v_3$ are consecutive vertices on $O_{\cal A}$ and none of the edges
$vv_1$, $vv_2$, $vv_3$ are claimed, then Avoider can first claim $vv_2$, and then one of the
edges $vv_1$, $vv_3$ in the following move, see Figure~\ref{f:ext}.

We observe that the number of bad vertices can never exceed five.
Indeed, if there were six bad vertices at some point of the game,
then there would be at least $6\lceil \frac{k}{3}\rceil>2k-4$ Enforcer's
edges.
That is more than the total number of edges played by Enforcer until that
point.

Let $k$ be the order of the current outerplanar graph $\cal A$. While $k< \lceil n/4\rceil$,
Avoider always attaches a good vertex of the highest Enforcer's degree. The setup behind this
process is similar to the one in the so-called Box Game.  Let $m$, $s$ and $\ell$ be positive
integers. In the Box Game, in each of the moves, the first player claims one element of the
board, and the second player claims $m$ elements of the board. The goal of the first player is
to claim one of $\ell$ disjoint winning sets of size $s$, and his opponent wants to prevent
him from doing that. We will make use of the following result.

\begin{theorem}[Chv\'atal and Erd\H{o}s $\cite{CE}$] \label{box}
The first player can win the Box Game when $s< m \log \ell$.
\end{theorem}

Avoider's strategy in our game keeps Enforcer's maximum degree over all good vertices minimal.
Theorem~\ref{box} with $m=4$ guarantees that this degree can never exceed $4\log n$.

Suppose that $k=\lceil n/4 \rceil$ and there are five bad vertices,
$b_1,\dots,b_5$.
We show how Avoider can reduce the number of bad
vertices to four in the two moves that follow.
As we have already seen, a bad vertex $v$ is adjacent in Enforcer's graph to
at least one of every three consecutive vertices on $O_{\cal A}$.
These vertices of $O_{\cal A}$
subdivide the edge set of $O_{\cal A}$ into paths of length at most
three, and we will refer to these as \emph{blocks}.
For every $i\in\{1,\dots,5\}$ and every edge $e$ on $O_{\cal A}$, we
define $f_i(e)$ as the set of edges in the block to which $e$
belongs, in the mentioned subdivision by $b_i$.
As we have seen, $|f_i(e)|$ is always either 1, 2 or 3.

The number of edges claimed by Enforcer between $b_i$ and
$V(O_{\cal A})$ is exactly $\sum_{e\in O_{\cal A}}
\frac{1}{|f_i(e)|}$, for $i\in\{1,\dots,5\}$.
The total number of edges claimed by Enforcer is not less than
$$
\sum_{i=1}^5 \sum_{e\in O_{\cal A}} \frac{1}{|f_i(e)|}.
$$
On the other hand, we know that Enforcer played at most $2k-4$ moves in
total, that is the number of edges in a $k$-degenerate graph, and hence,
$$
\sum_{e\in O_{\cal A}} \sum_{i=1}^5 \frac{1}{|f_i(e)|}\leq 2k-4.
$$
Therefore, there exists an edge $e\in O_{\cal A}$ such that $\sum_{i=1}^5 \frac{1}{|f_i(e)|} < 2$.
This can only happen if at least four of $|f_i(e)|$, $i=1,\dots,5$, say the first four, are
equal to 3.
Therefore, there has to be an edge $f$ on
$O_{\cal A}$ incident to $e$ such that $\{e,f\}$ belong to two of the blocks $f_i(e)$,
$i=1,\dots,5$, say, $f_1(e)$ and $f_2(e)$.

By $w_1,w_2,w_3$ we denote the three consecutive vertices on
$O_{\cal A}$ with $e=w_1w_2$, $f=w_2w_3$. Since $k=\lceil n/4
\rceil$, there still exists an isolated vertex $u$ in Enforcer's
graph. In the following move, Avoider claims the edge $uw_2$.

If Enforcer does not claim $uw_1$ in his response, Avoider claims it
immediately.
The vertex $u$ is also on $O_{\cal A}$ now, and four blocks
$f_i(e)$, $i=1,\dots,4$, are extended to size four in this way.
Only two of them can be subdivided by the last two Enforcer's moves.
Hence, some $b_i$ is not bad any more.

On the other hand, if Enforcer claims $uw_1$ in his response, then
Avoider claims $uw_3$, and similarly as before, blocks $f_1(e)$ and
$f_2(e)$ are extended to size four.
Enforcer can subdivide at most one of them in his following move, and
the bad vertex corresponding to the other block is not bad any more.

Therefore, after this process there are at most four bad vertices. As
long as $k<n-4$, Avoider keeps attaching good vertices with highest
Enforcer's degree to ${\cal A} \cup \{m\}$.
Since Enforcer's maximum degree over all good vertices is at most $4\log n$,
no other vertex can ever become bad.
Finally, when there are only bad vertices left, they are isolated in Avoider's graph.
Therefore, Avoider can play at least four more moves without creating a non-outerplanar
graph, and the total number of Avoider's moves is at least
$2(n-4)-4+4=2n-8$. {\hfill $\Box$ \medskip\\}


\subsection{Keeping the graph diamond-free}

\begin{figure}[ht]
 \begin{center}
  \includegraphics[scale=0.5]{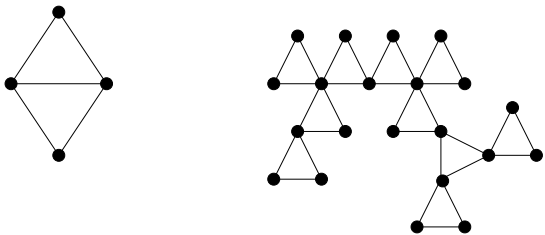}
 \end{center}
  \caption{A diamond, and a maximal diamond-free graph on 21 vertices}
\label{f:diam}
\end{figure}

If $G$ is a diamond-free graph on $n$ vertices, then the number of edges in $G$ is at most
$d(n)=\lceil\frac{3n-5}{2}\rceil$. We show that Avoider can survive in the game for nearly
that many moves. In Figure~\ref{f:diam} we see the diamond graph, and an example of a
diamond-free graph maximizing the number of edges.

{\it Proof of Theorem~\ref{th:df}.} The upper bound follows from the simple pairing argument
at the beginning of the proof of Theorem~\ref{th:op}.

For the lower bound, we give an explicit strategy for Avoider that enables him to play for
$d(n)-3$ moves. Before performing a detailed analysis, let us first sketch Avoider's strategy.
The game is divided into two phases. In the first phase, Avoider fixes two arbitrary vertices
$c_1$ and $c_2$ and connects them in his first move. Then, by using a pairing strategy, he
creates a spanning tree, consisting of two stars centered at $c_1$ and $c_2$, and the edge
$c_1 c_2$.
While doing that, Avoider pays attention to certain edge densities in Enforcer's graph,
preparing the ground for the second phase. In the second phase, Avoider claims a large
matching on the leaves of each of the stars. In this way, he forms a bunch of
edge-disjoint triangles along with a bridge and possibly some hanging edges, that is, a
diamond-free graph.

Next, we describe the first phase in detail.
Let $c_1$ and $c_2$ be two vertices, fixed before the game starts.
Avoider creates two disjoint stars centered in $c_1$ and $c_2$.
Throughout this phase, we denote the set of vertices adjacent to $c_i$ in Avoider's graph
by $L_i$, for $i=1,2$.
The set of vertices that are isolated in Avoider's graph is denoted by $R$.

\begin{figure}[ht]
 \begin{center}
  \includegraphics[scale=1.2]{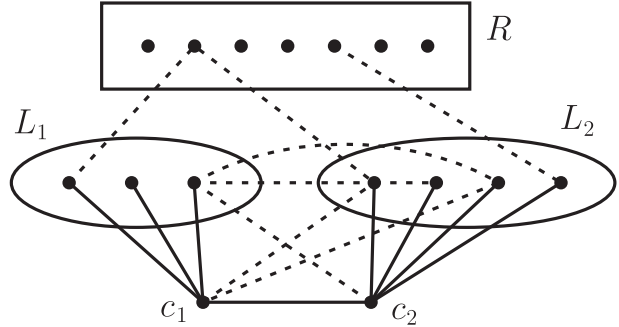}
 \end{center}
  \caption{A possible arrangement during the first phase -- solid lines represent Avoider's edges,
  and dashed lines represent Enforcer's edges}
\label{f:setup}
\end{figure}

We list the rules for Avoider's strategy during the first phase. In the first move, he claims
the edge $c_1 c_2$. The rest of the rules follow. One possible arrangement of edges played is
shown in Figure~\ref{f:setup}.

\begin{itemize}
 \item
Whenever Enforcer claims an edge $xc_i$, for some $i\in\{1,2\}$,
Avoider responds by claiming the edge $xc_{3-i}$,
 \item
If Enforcer claims an edge $uv$, where $u\in L_i$, for some
$i\in\{1,2\}$, and $v\in R$, then Avoider responds with $vc_{3-i}$,
 \item
If Enforcer claims an edge $uv$, where $u,v\in L_i$, for some
$i\in\{1,2\}$, then Avoider responds with $wc_{i}$, for arbitrary
$w\in R$,
 \item
If Enforcer claims an edge $uv$, with $u,v\in R$, then Avoider
responds by claiming $c_iu$, where $i$ is arbitrary,
\item
If Enforcer claims an edge between $L_1$ and $L_2$, then Avoider
responds by claiming $c_iu$, where $u$ is any vertex from $R$, and
$i$ is arbitrary.
\end{itemize}
The first phase ends as soon as Avoider plays a move after which
$V=L_1\cup L_2 \cup \{c_1,c_2\}$ holds.

Let $E_{\cal E}(X)$ denote the set of edges in Enforcer's graph, induced
by $X$.
We define the following density measure,
\begin{equation}
\varrho_i=\frac{|E_{\cal E}(L_i)|+|E_{\cal
E}(L_i,R)|}{\max\{|V(L_i)|, 1\}} \text{, for }i=1,2. \label{e:dens}
\end{equation}
We prove that throughout the first phase, after
each of his moves, Avoider keeps both $\varrho_1$ and $\varrho_2$ to be at most 1.
Indeed, the densities from (\ref{e:dens}) are certainly less than 1
after the very first move of Avoider. Next, let us look at a move of
Enforcer, and the corresponding move of Avoider. Checking through
all of the rules in Avoider's strategy we see that the densities
either remain unchanged, or 1 is added to both the numerator and the
denominator in (\ref{e:dens}).
Hence, neither of the densities can exceed 1.

We proceed to the second phase, in which Enforcer plays the first move.
As we have already mentioned, Avoider's goal is to build a
large matching on both $L_1$ and $L_2$.
Throughout this phase, as soon as Avoider claims an edge $v_1v_2\in L_i$, for some
$i\in\{1,2\}$, we remove both $v_1$ and $v_2$ from $L_i$. The
set of rules for Avoider's strategy in this phase can now be
described.
\begin{itemize}
 \item If Enforcer claims an edge in $L_i$, for some
$i\in\{1,2\}$, Avoider responds by claiming an edge also in $L_i$.
Otherwise, Avoider claims an edge in any of the sets $L_i$, $i\in\{1,2\}$.
 \item Whenever Avoider wants to respond by playing in $L_i$,
 $i\in\{1,2\}$, we distinguish two cases:
\begin{enumerate}
 \item If there is an unclaimed edge in $L_i$ that is adjacent to
 a vertex $m$ with maximum Enforcer's
 degree in $L_i$, Avoider claims it.
 \item If there is no unclaimed edge in $L_i$ that is adjacent to
 the vertex $m$ with maximum Enforcer's
 degree in $L_i$, Avoider removes $m$ from $L_i$, and then claims an
 edge following again this set of rules.
 \end{enumerate}
\end{itemize}

Following these rules Avoider keeps $\varrho_i\le 1$, where $R$ in (\ref{e:dens})
is now the empty set, for $i=1,2$.
If the above condition in 2. is satisfied, knowing that $\varrho_i\le 1$,
the new set $L_i$ induces at most one Enforcer edge, $e$ say.
Avoider's reply is an edge incident to $e$, and therefore Enforcer's graph $\cal E$ becomes
empty on $L_i$, and remains empty after every of the following moves
of Avoider during phase two.
That is, case 2. above can happen at most once for each $L_i$, and whenever it happens Avoider can reach
a near-perfect matching in that $L_i$.

If case 2. does not occur for $L_i$, $i\in\{1,2\}$, Avoider can
follow the algorithm until $|L_i|<4$. If $|L_i|\le 2$, then Avoider
has reached a near-perfect matching.
The only case when Avoider is possibly stuck is $|L_i|=3$, if after Enforcer's move,
$L_i$ spans a triangle of Enforcer's edges.
We conclude that the total number of vertices in $L_1\cup L_2$ that are unsaturated by the two
matchings is at most six.

As we have already mentioned, phase two is finished, when the matchings of Avoider
can not be further extended.
By simply counting the edges played, the lower bound from the theorem readily
follows. {\hfill $\Box$ \medskip\\}


\subsection{Keeping the graph $k$-degenerate}

Before we prove the theorem, notice an alternative way of defining degeneracy:
a graph $G$ is $k$-degenerate if and only if there is a total ordering of
$V(G)$ such that any vertex has at most $k$ preceding neighbors in that ordering.

{\it Proof of Theorem~\ref{th:kd}.} We exhibit a strategy for Avoider to claim the edges of a
maximal $k$-degenerate graph on $n$ vertices in his first $e(n)$ moves. We split the game into
two phases.

In the first phase, Avoider wants to create a maximal $k$-degenerate graph on significantly
less than $n$ vertices. In the second phase, he gradually attaches all the remaining vertices
to that graph.

Let us now describe both phases in detail. The first phase is subdivided into $k$ subphases.
In the beginning of the first subphase, Avoider picks a vertex $v_1$, and he repeatedly claims
edges adjacent to $v_1$ until he has claimed $3^{3k}$ edges. By $V_1$ we denote the set of
vertices adjacent to $v_1$ in Avoider's graph at this point. For $2\leq i \leq k$, in the
beginning of the $i$-th subphase, Avoider chooses the vertex $v_i\in V_{i-1}$ of minimal
degree in Enforcer's graph induced on $V_{i-1}$, and connects $v_i$ to some $3^{3k-i+1}$
vertices of $V_{i-1}$. We denote the set of those vertices by $V_i$. It remains to show that
this can be done, i.e., before the $i$-th subphase there are at least $2\cdot 3^{3k-i+1}$
vertices in $V_{i-1}$ such that edges between them and $v_i$ are not claimed. This can be seen
as follows. The total number of moves played in the first $i-1$ subphases is $\sum_{j=1}^{i-1}
3^{3k-t+1} \leq 3^{3k+1}/2$, and the minimum degree in the Enforcer's graph taken over all
vertices in $V_{i-1}$ is not greater than $3^{3k+1}/|V_{i-1}|=3^{i-1}$, implying our claim.

After the end of the first phase, Avoider's graph induced on the set $R=V_k \cup
\{v_1,\dots,v_k\}$ is a maximal $k$-degenerate graph. During the second phase, the vertices
from $V(G)\setminus R$ will be gradually attached to that graph using a pairing strategy. Note
that Avoider has claimed some edges between $R$ and $V\setminus R$ already in the first phase.

For every vertex $x\in V(G)\setminus R$, Avoider's first hope is to claim $k$ edges between
$x$ and $R$, including the edges played in the first phase. To check if that can be done using
a pairing strategy, to each $x$ we assign the following number,

$$f(x):=deg_{\cal A} (x,R) + \frac{1}{2}(|R|- deg_{\cal E} (x,R) - deg_{\cal A} (x,R)).$$

Here, $deg_{\cal A} (x,R)$ and $deg_{\cal E} (x,R)$ stand for the numbers of edges between $x$
and $R$ claimed by Avoider and Enforcer, respectively. By $D$ we denote all vertices in
$V(G)\setminus R$ with $f(x)\geq k$, and let $F:=\left(V(G)\setminus R\right) \setminus D$.
Since the total number of edges claimed in the first phase is less than $3^{3k+1}$, we know
that $|F|\leq 3^{3k+1} <n/2$.

Now, for every vertex $v\in D$, Avoider will use a simple pairing strategy to claim $k$ edges
between $v$ and $R$, also counting the edges he has already claimed in the first phase. To do
that, he considers $2(k-deg_{\cal A} (x,R))$ unclaimed edges between $v$ and $R$, and pairs
them up arbitrarily.

For every vertex $v\in F$, Avoider aims at connecting it to a larger set, $R\cup D$. He will
again use a simple pairing strategy to claim $k$ edges between $v$ and $R\cup D$. To do that,
he considers $2(k-deg_{\cal A} (x,R))$ unclaimed edges and pairs them up arbitrarily.

Avoider's strategy for the second phase is the following. Whenever Enforcer claims one of the
paired edges, A\-void\-er immediately responds by claiming the other one. If Enforcer claims
an edge that does not belong to a pair, then Avoider claims an edge in an arbitrary pair, and
removes that pair for the rest of the game. As long as Avoider proceeds like this, he will not
lose. Indeed, looking at the alternative definition of $k$-degeneracy presented in the
beginning of this section, we see that any total ordering $\prec$ in which $\{v_1,\dots,v_k\}
\prec V_k \prec D \prec F$ verifies that Avoider's graph is $k$-degenerate. When all the pairs
are removed he has already claimed a maximal $k$-degenerate graph on $n$ vertices. {\hfill
$\Box$ \medskip\\}

\section{Concluding remarks and open problems}
\label{sec::openprobs}

Looking at the Avoider-Enforcer diamond-free game, and the games of non-planarity and
non-outerplanarity, we could observe a pattern regarding how long the game lasts. Namely, the
number of moves Avoider can survive in those games are all just an additive constant away from
the upper bound in (\ref{e:ob}). We are curious whether this pattern extends to a larger class
of forbidden graphs.

\begin{question}
Let $H$ be a fixed graph, and let ${\cal F}^H_n$ be the set of
subgraphs of $K_n$ that contain an $H$-minor. Is it true that
$$
\tau_E({\cal F}^H_n) = \text{\rm ex}({\cal F}^H_n) +\BigO (1)?
$$
\end{question}

Even though our main goal was to prove Theorem~\ref{th:kd} for constant values of $k$, it
turned out that our proof readily holds for all $k=\littleo (\log n)$. We did not make
particular efforts to analyze the same problem for larger values of $k$. Still, we think that
it would be interesting to find out for how large $k$, in terms of $n$, the statement of
Theorem~\ref{th:kd} still holds.

\begin{question}
How large can $k=k(n)$ be, so that $\tau_E({\cal D}_n^k)=e(n)+1$ still holds?
\end{question}


\end{document}